\documentclass[12pt]{article}
\usepackage{amsmath,amsthm,latexsym}
\input{epsf.sty}
\makeatletter
\newfont{\footsc}{cmcsc10 at 8truept}
\newfont{\footbf}{cmbx10 at 8truept}
\newfont{\footrm}{cmr10 at 10truept}
\renewcommand{\ps@plain}{%
\renewcommand{\@oddfoot}{\footsc the electronic journal of combinatorics
  {\footbf 10} (2003), \#R00\hfil\footrm\thepage}}
\makeatother
\newcommand{\eqdef}{\, =\kern -12.7pt\raise 6pt\hbox{{\tiny\textrm{def}}}\,\,}

\newtheorem{theorem}{Theorem}
\numberwithin{equation}{section}
\begin{document}
\title{Longest increasing subsequences\\ in pattern-restricted permutations}
\author{
Emeric Deutsch\\
\small Polytechnic University\\
\small Brooklyn, NY 11201\\
\small \texttt deutsch@duke.poly.edu
\and
A. J. Hildebrand\\
\small University of Illinois\\
\small Urbana, IL 61801\\
\small \texttt{hildebr@math.uiuc.edu}
\and
Herbert S. Wilf\\
\small University of Pennsylvania\\
\small Philadelphia, PA 19104-6395\\
\small {\texttt wilf@math.upenn.edu}
}
\maketitle

\begin{center}
\date{\small Submitted: April 8, 2003;  Accepted: April 21, 2003.}\\
{\small MR Subject Classifications: 05A16, 05A05}
\end{center}

\begin{abstract}
Inspired by the results of Baik, Deift and Johansson on the limiting distribution of the lengths of the longest increasing subsequences in random permutations, we find those limiting distributions for pattern-restricted permutations in which the pattern is any one of the six patterns of length 3. We show that the (132)-avoiding case is identical to the distribution of heights of ordered trees, and that the (321)-avoiding case has interesting connections with a well known theorem of Erd\H os-Szekeres.
\end{abstract}

\section{Introduction}
A great deal of spectacular work has been done on determining the distribution of the length of the longest increasing subsequence, LIS($p$), in a random permutation $p$. This work culminated in the paper of Baik, Deift and Johansson \cite{bdj} which found the limiting distribution function. 

Historically, the \textit{average} length, $\mu(n)$, of the LIS had been found  to be $\ge \sqrt{n}/2$ \cite{es} and then  $\ge 2\sqrt{n}$ \cite{ls}. The correct order of magnitude was established in \cite{vk}, where it was shown to be $\sim 2\sqrt{n}$. Next the standard deviation $\sigma(n)$ had been estimated \cite{amo} to be approximately $n^{1/6}$ but was not proved to be so. Finally, in \cite{bdj}, the complete limiting distribution function of the normalized random variable $(\mathrm{LIS}-\mu(n))/n^{1/6}$ was determined.

We can ask the same questions in some given subset of the set of all permutations. For example, what is the distribution of the length of the longest increasing subsequence in a randomly chosen permutation from the set of those that avoid the pattern\footnote{A permutation $p$ avoids (132) (resp. (231), (321)) if there do not exist $1\le i<j<k\le n$ s.t. $p(i)<p(k)<p(j)$ (resp. $p(k)<p(i)<p(j)$, $p(i)>p(j)>p(k)$).} $(231)$, or the pattern $(321)$? The tools are available for the study of such problems. A paper of Reifegerste \cite{re} finds a simple expression for the frequencies of different LIS lengths in the $(231)$-avoiding case. For $(132)$-avoiding permutations it turns out that the distributions are identical with known distribution functions of ordered trees by height, so we can simply quote the known results. Finally, for $(321)$-avoiding permutations the problem is related to a celebrated theorem of Erd\H os-Szekeres. In this case we exhibit the relevant statistics and find explicitly the limiting distribution function.

We remark that for all six of the patterns of length 3 it is true that the number of permutations of $n$ letters that avoid the pattern is the Catalan number $C_n={2n\choose n}/(n+1)$. So for each of the six patterns that we will study here, we will be asking what fraction of the roughly $c4^n/n^{3/2}$ permutations that avoid that pattern have a LIS of some given length.

 Our main results are the following.
\begin{theorem}
In the class of \,$(231)$-avoiding permutations of $n$ letters, and in the class of $(312)$-avoiding permutations of $n$ letters, the length of the longest increasing subsequence has mean $(n+1)/2$,  and standard deviation $\sim \sqrt{n}/2$. Moreover, the random variable $(LIS(p)-(n+1)/2)/(\sqrt{n}/2)$ has asymptotically the standard normal distribution.
\end{theorem}
\begin{theorem}
In the class of \,$(132)$-avoiding permutations of $n$ letters, and in the class of $(213)$-avoiding permutations of $n$ letters,  the length of the longest increasing subsequence has mean $= \sqrt{\pi n}+O(n^{\frac{1}{4}})$ and standard deviation 
\begin{equation}
\label{eq:stdv}
=c_1\sqrt{n}+O(n^{\frac14})\qquad\left(c_1=\sqrt{\pi\left(\frac{\pi}{3}-1\right)}=0.38506..\right).
\end{equation}
Moreover, the normalized random variable
\[X_n(p)\eqdef \frac{\mathrm{LIS}(p)-\sqrt{\pi n}}{\sqrt{n}},\]
defined for $n$-permutations $p$ that avoid the pattern (132), satisfies
\begin{equation}
\label{eq:limdist}
\lim_{n\to\infty}\mathrm{Prob}(X_n(p)\le \theta)=\sum_{t=-\infty}^{\infty}(1-2t^2(\theta+\sqrt{\pi})^2)e^{-(\theta+\sqrt{\pi})^2t^2}\qquad (\theta>-\sqrt{\pi}).
\end{equation}
\end{theorem}
\begin{theorem}
For permutations $p$ in the class of $(321)$-avoiding permutations of $n$ letters, define the random variable
\[X_n(p)\eqdef \frac{\mathrm{LIS}(p)-\frac n2}{\sqrt{n}}.\]
Then we have
\begin{eqnarray*}
\lim_{n\to\infty}\mathrm{Prob}(X_n(p)\le \theta)&=&\frac{2}{\sqrt{\pi}}\int_0^{4\theta^2}u^{1/2}e^{-u}du\\
&=&\frac{\Gamma(3/2,4\theta^2)}{\Gamma(3/2)},
\end{eqnarray*}
where $\Gamma(z,w)$ is the incomplete Gamma function
\[\Gamma(z,w)=\int_0^{w}u^{z-1}e^{-u}du.\]
\end{theorem}

\section{Avoiding the pattern $(231)$}

One situation where it is simple to deduce the limiting distribution is that of $(231)$-avoiding permutations.  This is because Reifegerste \cite{re} has shown that the number of $(231)$-avoiding permutations whose longest increasing subsequence has length exactly $k$ is
\begin{equation}
\label{eq:reif}
e(n,k)=\frac{1}{n}{n\choose k}{n\choose k-1}.
\end{equation}
The number of all $(231)$-avoiding permutations of $n$ letters is $C_n$, the Catalan number.
From this it is easy to check that the mean length of the LIS in this family is $(n+1)/2$ and the standard deviation is $\sqrt{n}/2+O(1)$. Finally it is a straightforward exercise in Stirling's formula to find that
\[\lim_{n\to\infty}\sqrt{n}\frac{e\left(n,\frac{n+1}{2}+\theta\frac{\sqrt{n}}{2}\right)}{C_n}=\frac{2}{\sqrt{\pi}}e^{-\theta^2}.\]
This says that for $(231)$-avoiding permutations $p$, the random variable
\[X(p)=\frac{LIS(p)-\frac{n+1}{2}}{\frac{1}{2}\sqrt{n}}\]
has asymptotically the standard normal distribution. 

It seems noteworthy that in this case the mean is quite large, i.e., typically such a permutation contains an increasing subsequence whose length is about half that of the permutation itself. 

It also is noteworthy that (\ref{eq:reif}) shows that the number of $(231)$-avoiding permutations of $n$ letters whose longest increasing subsequence has length $k$ is the same as the number whose longest increasing subsequence has length $n-k+1$. It is arresting that the number of $(231)$-avoiding permutations of, say, 1000 letters, whose LIS has length 980 is equal to the number whose LIS has length 21.

\section{Avoiding the pattern $(132)$}
The case of $(132)$-avoiding permutations is handled by the following sequence of observations.
\begin{itemize}
\item The number $f(n,k)$ of $n$-permutations that avoid $(132)$ and whose longest increasing subsequence has length $< k$ is equal to the number of ordered trees of $n$ edges whose height is $< k$. This follows from a bijection of Jani and Rieper \cite{jr}. 
\item It is known from the theory of ordered trees \cite{bkr} that this number is 
\begin{equation}
\label{eq:ffexpl2}
f(n,k)=2\sum_{t=-\infty}^{\infty}\left({2n\choose n+t(k+1)}-\frac{1}{4}{2n+2\choose n+1+t(k+1)}\right)\quad (k=1,2,3,\dots).
\end{equation}
\item To see that the normalized random variable
\[X(p_n)=\frac{\mathrm{LIS}(p_n)-\sqrt{\pi n}}{\sqrt{n}},\]
defined on $(132)$-avoiding permutations $p_n$ of $n$ letters, has a nontrivial limiting distribution, as $n\to\infty$, we consider
\begin{eqnarray*}
\mathrm{Prob}(X(p_n)< \theta)&=&\mathrm{Prob}(\mathrm{LIS}(p_n)\le (\theta+\sqrt{\pi}\,)\sqrt{n})\\
&=&\frac{f(n,\sqrt{n}(\theta+\sqrt{\pi}\,))}{C_n},
\end{eqnarray*}
where
\[C_n=\frac{1}{n+1}{2n\choose n}\]
is the  $n$th Catalan number. But by (\ref{eq:ffexpl2}) we see that
\[f(n,\sqrt{n}(\theta+\sqrt{\pi}\,))=\frac{1}{C_n}\sum_{t=-\infty}^{\infty}\phi(n,t),\]
in which the general $t$-th term of the sum is
\[\phi(n,t)=\frac{n+1}{{2n\choose n}}{2n\choose n+tc\sqrt{n}}\left(2-\frac{(n+1)(2n+1)}{(n+1)^2-t^2c^2n}\right).\qquad(c=\theta+\sqrt{\pi})\]
It is an easy exercise in Stirling's formula to see that
\[\lim_{n\to\infty}\phi(n,t)=(1-2c^2t^2)e^{-c^2t^2}\]
and the proof is complete. 
\end{itemize}

The limiting cumulative distribution function $F(\theta)$, in (\ref{eq:limdist}), is related to the Jacobi theta function
\[\Theta(x)=\sum_{t=-\infty}^{\infty}e^{-t^2\pi x}.\]
In fact, $F(\theta)=\Theta(u)+2u\Theta'(u)$, where $u=(1+\theta/\sqrt{\pi})^2$. It is well known \cite{fo} that this is indeed the limiting distribution, from the theory of ordered trees, but the derivation is simple enough that we have included it above. Indeed \cite{fo} shows that this same distribution function, originally considered by R\'enyi and Szekeres \cite{rs}, is the limiting distribution function for many classes of trees.

We remark that although the Jani-Rieper bijection \cite{jr} is the most explicit mapping between $132$-avoiding permutations of $n$ letters whose LIS has length $k$ and ordered trees of $n$ edges whose height is $k$, there are other, somewhat less explicit, bijections in the literature. Indeed, in Krattenthaler \cite{kr}, we find an explicit bijection between such permutations and Dyck paths. This mapping has the property that the length of the LIS corresponds to the height of the image path. In view of well known bijections between Dyck paths by height and ordered trees by height, we have, by composition, the mapping that we need here. Likewise in Chow, West \cite{ChWe} there are remarks in section 5 that suggest, less explicitly, that such bijections exist. Other, equivalent bijections with Dyck paths have been given by Fulmek \cite{fu}, Reifegerste \cite{re}, and Bandlow and Killpatrick \cite{bk}.

\section{The $(321)$ case and the Erd\H os-Szekeres theorem}

Finally we ask for the limiting distribution function for the length of the longest increasing subsequence in the class $S_n(321)$ of $(321)$-avoiding permutations of $n$ letters. Let $f(n,k)$ denote the number of such permutations whose longest increasing subsequence has length $\le k$. A famous theorem of Erd\H os-Szekeres states that every permutation of $(r-1)(s-1)+1$ or more letters contains either a decreasing subsequence of length $r$ or an increasing subsequence of length $s$. If we take $r=3$ we see that every permutation in the class $S_{2k+1}(321)$ has an increasing subsequence of length $k+1$, i.e. $f(n,k)=0$ for all $n\ge 2k+1$.

Even though, for each fixed $k$, $f(n,k)$ vanishes for all large enough $n$, the question of the limiting distribution function of the LIS remains.
\begin{theorem}
\label{th:321}
The number $a(n,k)$ of $(321)$-avoiding permutations having longest increasing subsequence of length $= k$ is given by
\begin{equation}\label{eq:ankform}
        a(n,k)= \left(\frac{2k-n+1}{n+1}{n+1\choose k+1}\right)^2
\end{equation}
for $\lfloor{(n+1)/2}\rfloor \le  k \le  n,$ and vanishes otherwise.
\end{theorem}

These permutations correspond 1-1, under the Schensted insertion algorithm (see, e.g., \cite{bs}), to pairs 
of Young tableaux with $n$ cells, at most two rows, and $k$ columns, so the number of
them is $B^2$, where $B$ is the number of such tableaux. A brief computation
with the hook formula shows that \[B=\frac{2k-n+1}{n+1}{n+1\choose n-k},\]
completing the proof. $\Box$

It follows that the number $f(n,k)$ of permutations in $S_n(321)$ whose LIS has length $\le k$ is
\begin{equation}
f(n,k)=\sum_{\lfloor{(n+1)/2}\rfloor \le  j \le  k}\left(\frac{2j-n+1}{n+1}{n+1\choose j+1}\right)^2.
\end{equation}
Now fix $t>0$. The probability that a permutation of $n$ letters has LIS of length at most $n/2+t\sqrt{n}$ is \[\frac{f(n,n/2+t\sqrt{n})}{\frac{1}{n+1}{2n\choose n}}.\]
Hence the limiting probability distribution function that we seek is
\begin{eqnarray*}
F(t)&=&\lim_{n\to\infty}\frac{1}{\frac{1}{n+1}{2n\choose n}}\sum_{\lfloor{(n+1)/2}\rfloor \le  j \le  n/2+t\sqrt{n}}\left(\frac{2j-n+1}{n+1}{n+1\choose j+1}\right)^2\\
&=&\lim_{n\to\infty}\frac{4}{(n+1){2n\choose n}}\sum_{0\le k\le t\sqrt{n}+1/2}k^2{n+1\choose (n+1)/2+k}^2.
\end{eqnarray*}
By an easy application of Stirling's formula the relation
\[\frac{{2m\choose m-k}}{{2m\choose m}}\sim e^{-k^2/m}\]
holds uniformly for  all $|k|< t\sqrt{m}$, and we have
\[F(t)=\frac{32}{\sqrt{\pi}}\lim_{n\to\infty}\frac{1}{\sqrt{n}}\sum_{0\le k\le t\sqrt{n}}\frac{k^2}{n}e^{-4k^2/n}.\]
The quantity whose limit is being taken is a Riemann sum that approximates the integral
\[\int_0^tx^2e^{-4x^2}\]
with grid size $1/\sqrt{n}$, and so the desired limit is
\begin{eqnarray*}
F(t)&=&\frac{32}{\sqrt{\pi}}\int_0^tx^2e^{-4x^2}dx\\
&=&\frac{2}{\sqrt{\pi}}\int_0^{4t^2}u^{1/2}e^{-u}du\\
&=&\frac{\Gamma(3/2,4t^2)}{\Gamma(3/2)}.
\end{eqnarray*}

The limiting distribution $F(t)$ is well-known in statistics: Except for a linear change of variables, it is a $\chi^2$ distribution
with three degrees of freedom, which arises as the distribution of the length of a $3$-vector whose entries are independent standard normal
random variables (see, for example, \cite{as}).  It would be interesting to find a heuristic explaining this connection. 

\section{Other patterns of three letters}
We have so far dealt with the statistics of the length of the longest increasing subsequence in permutations that avoid any one of the patterns (132), (231), or (321). There remain three other patterns of three letters, viz. (123), (312), and (213).

The (123) case is trivial. If an $n$-permutation avoids (123) then its LIS has length exactly 2 except for the reversal of the identity permutation, whose LIS has length 1.

The (312)  case is identical with the (231). Indeed, the length of the LIS of a (312)-avoiding permutation is equal to  the length of the longest decreasing subsequence of its reversal, which is (213)-avoiding, and that is equal to the length of the LIS of the complement of its reversal, which is (231)-avoiding, as claimed.
                
Similarly, by reversal and complementation, the (213) case is identical with the (132).

\vspace{.3in}

We thank Christian Krattenthaler, Brendan McKay and Julian West for informative discussions related to this paper.

\end{document}